\documentclass{article}
\usepackage{amsfonts,amssymb,amsmath}

\newenvironment{proof}[1][Proof]{\textbf{#1.} }{\hfill
  \rule{0.5em}{0.5em} \medskip}

\newcommand{\cl}[1]{\overline{#1}}

\newtheorem{theorem}{Theorem}

\newtheorem{lemma}[theorem]{Lemma}

\newtheorem{cor}[theorem]{Corollary}

\newtheorem{prop}[theorem]{Proposition}

\newtheorem{ques}[theorem]{Question}

\newcommand{\R}{\mathbb{R}}
\newcommand{\N}{\mathbb{N}}

\newcommand{\plim}{\underleftarrow{\lim}}

\newcommand{\bfam}{{\Phi}}

\newcommand{\basic}{\text{basic}\,}

\newcommand{\cof}{\mathop{cof}} 

\newcommand{\ctm}{\mathbf{c}}

\begin{document}

\title{Minimal Size of Basic Families}

\author{Ziqin Feng and Paul Gartside}

\date{July  2009}

\maketitle

\begin{abstract} A family $\bfam$ of continuous real-valued
functions on a space $X$ is said to be {\sl basic} if  every $f \in C(X)$ can be
represented $f = \sum_{i=1}^n g_i \circ  \phi_i$ for some $\phi_i \in \bfam$
and $g_i \in C(\R)$ ($i=1, \ldots , n$). Define $\basic (X) = \min \{ |\bfam| : \bfam$ is a basic family for $X\}$. If $X$ is separable metrizable $X$ then either $X$ is locally compact and finite dimensional, and $\basic (X) < \aleph_0$, or  $\basic (X) = \mathfrak{c}$. If $K$ is compact and either $w(K)$ (the minimal size of a basis for $K$) has uncountable cofinality or $K$ has a discrete subset $D$ with $|D|=w(K)$  then either $K$ is finite dimensional, and $\basic (K) = \cof ([w(K)]^{\aleph_0} , \subseteq)$, or $\basic (K) = |C(K)|=w(K)^{\aleph_0}$.
\end{abstract}

\section{Introduction}

The 13th Problem of Hilbert's celebrated list \cite{H} asks whether every continuous real valued function of three variables can be
written as a superposition (i.e. composition) of continuous functions
of two variables. Hilbert conjectured that the answer was no, but in
1957 Kolmogorov, building on previous work of himself  and
Arnold,
proved a remarkable  result:
 every continuous real valued function of
$n$--variables from a closed and bounded interval can be expressed as a
superposition of functions of just one variable,  and addition.
\begin{theorem}[Kolmogorov Superposition, \cite{Kol}]\label{KST}
For a fixed $n \ge 2$, there are $n(2n+1)$ continuous maps $\psi_{pq}:[0,1] \to
\mathbb{R}$  such that every continuous $f:[0,1]^n \to \mathbb{R}$ can be
written: 
\[ f(\mathbf{x}) = \sum_{q=1}^{2n+1} (g_q \circ \phi_q) (\mathbf{x}) 
\quad \text{where } \phi_q(x_1, \ldots , x_n) = \sum_{p=1}^n
\psi_{pq} (x_p),\] 
and the $g_q:\mathbb{R} \to \mathbb{R}$ are continuous maps depending on $f$.
\end{theorem} 
Recently the authors have extended the Kolmogorov Superposition Theorem to functions of $n$ {\sl real} variables, \cite{FG2}. This gives a more complete solution of Hilbert's 13th Problem free from the restriction to bounded intervals which is unnatural in the context (solution functions of polynomials) that Hilbert placed his problem.

In this paper we focus on the functions $\phi_q$ in Kolmogorov's theorem. Even in the case when $n=2$, Theorem~\ref{KST} says something unexpected and insightful: there are just $5$ continuous functions, $\phi_1, \ldots , \phi_5$, on the unit square so that {\sl every} continuous function on the unit square can be obtained in a simple way from just these $5$ functions along with functions of one real variable. In other words, to understand $C([0,1]^2)$ it suffices to understand $C([0,1])$ and these $5$ functions. (Here and below, all topological spaces are Tychonoff, $C(X,Y)$ is the set of all continuous functions from $X$ to $Y$, and $C(X)=C(X,\R)$. We write $\aleph_\alpha$ for the $\alpha$th infinite cardinal and $\mathfrak{c}$ for $2^{{\aleph_0}}$.)

Following Sternfeld and others a family $\bfam$ of continuous real-valued
functions on a space $X$ is said to be {\bf basic} if  every $f \in C(X)$ can be
represented $f = \sum_{i=1}^n g_i \circ  \phi_i$ for some $\phi_i \in \bfam$
and $g_i \in C(\R)$ for $i=1, \ldots , n$.

In \cite{FG1}, building on work of Sternfeld \cite{St}, Ostrand \cite{Ost}, and others, the authors showed that a space $X$ has a finite basic family if and only if $X$ is locally compact, finite dimensional and separable metrizable (or equivalently, homeomorphic to  a closed subspace of Euclidean space). In this case, $\dim (X) \le n$ if and only if $\basic (X) \le 2n+1$. It might seem plausible that a space $X$ has a countable basic family precisely when $X$ is suitably `nice' and countable dimensional, but this is not the case. The result from \cite{FG1} says that if a space has a countable basic family, then in fact it has a finite basic family.

These results help motivate the following definition of a new cardinal invariant of topological spaces: $\basic (X) = \min \{ |\bfam| : \bfam$ is a basic family for $X\}$. Natural questions arise:  what are the possible values for $\basic (X)$? can we calculate, or at least bound,  $\basic(X)$ using other
cardinal invariants of $X$, such as weight, $w(X)$, the minimal size of a basis for $X$? 

Since the natural map of $X$ into $\R^\bfam$ is an embedding when $\bfam$
 is a basic family a simple restriction on the size of basic families is: $ w(x) \le \basic (X) . {\aleph_0}  \le |C(X)|$. So further natural questions are: when is $\basic (X) \le w(X)$? when is $\basic (X) = |C(X)|$? is it possible to have $\basic(X)$ strictly between $w(X)$ and $|C(X)|$?

In this paper we consider these questions for {\sl separable metrizable} spaces and {\sl compact} spaces. Suppose first that $X$ is separable metrizable. Then from the above, either $\basic (X)$ is finite, and this happens if and only if $X$ is locally compact and finite dimensional, or $\aleph_1 \le \basic(X) \le \mathfrak{c} = |C(X)|$. Experience of other related cardinal invariants of separable metrizable spaces would suggest that $\basic (X)$ should be undetermined by the standard axioms of set theory (ZFC). For example $k(X)$, which is the minimal size of a cofinal family in the set of all compact subsets of $X$, is undetermined even when $X$ is the rationals or the irrationals. However (Theorem~\ref{sm}) $\basic (X)$ {\sl is} determined in ZFC for all separable metrizable $X$: 
\begin{itemize}
\item[either] $X$ is locally compact and finite dimensional, and $\basic (X) < \aleph_0$,
\item[or] $X$ is either infinite dimensional or not locally compact, and  $\basic (X) = \mathfrak{c}$.
\end{itemize}

This theme --- that $\basic (X)$ is remarkably absolute --- is continued when we consider compact spaces. Note that if $K$  is compact, then Stone \cite{Stone} has shown that $|C(K)| = w(K)^{{\aleph_0}}$.
Hence,  $\basic (K)$ lies between the weight of $K$ and the countable power of the weight. This leads to some intriguing connections with Shelah's {\sl Potential Cofinalities Theory} (PCF).

Let $\kappa$ be an uncountable cardinal.  Shelah observed that $\kappa^{\aleph_0} = \cof ([\kappa]^{\aleph_0}, \subseteq) \times | \mathbb{P}({\aleph_0})|$. (Here $\cof ([\kappa]^{\aleph_0}, \subseteq)$ is the minimal size of a cofinal set in the countably infinite subsets of $\kappa$ ordered by inclusion.)
If $\kappa$ has uncountable cofinality then $\cof ([\kappa]^{\aleph_0}, \subseteq) = \kappa$, and so $\kappa^{\aleph_0}$ is easily computed --- it is $\max (\kappa, \mathfrak{c})$.

However, if $\kappa$ has countable cofinality then Shelah has shown \cite{Sh} that interesting things happen. Whereas the value of $|\mathbb{P}({\aleph_0})|=\mathfrak{c}$ is almost entirely unconstrained by the axioms of set theory and can be made arbitrarily large, $\cof ([\kappa]^{\aleph_0}, \subseteq)$ seems to be almost absolute. For example $\aleph_{\omega} < \cof ([\aleph_{\omega}]^{\aleph_0}, \subseteq) < \aleph_{{\omega}_4}$, and making $\cof ([\aleph_{\omega}]^{\aleph_0}, \subseteq) > \aleph_{{\omega}+1}$ requires large cardinals.

We prove (Theorems~\ref{cpt_dichot1} and~\ref{cpt_dichot2}) that if $K$ is compact and either $w(K)$ has uncountable cofinality or $K$ is suitably `nice' then
\begin{itemize}
\item[either] $K$ is finite dimensional, and $\basic (K) = \cof ([w(K)]^{\aleph_0} , \subseteq)$,
\item[or] $K$ is infinite dimensional, and $\basic (K) = |C(K)|=w(K)^{\aleph_0}$.
\end{itemize}
This gives almost complete information on the possible values of $\basic (K)$ for compact $K$. These are teased out and examples given below.

It is also interesting to note that if $K$ is compact, finite dimensional, `nice' and of weight $\kappa$ (for example, $K=2^\kappa$), and if $\bfam$ is a basic family for $K$ of minimal size, then $C(K) \sim \bigcup_{n \in \N} \left(\bfam^n \times C(\R)^n \right)$ is a natural `topological realization' of the cardinal identity $\kappa^{\aleph_0} = \cof ([\kappa]^{\aleph_0}, \subseteq) \times | \mathbb{P}({\aleph_0})|$.

Finally we briefly discuss connections of the above results with Banach algebras.
Let $K$ be a compact space. Then $C(K)$ with the supremum norm is a Banach algebra. Sternfeld has observed that for any $\phi \in C(K)$ the set $L(\phi) = \{ g \circ \phi : g \in C(\R)\}$ is a closed subring of $C(K)$ containing the constants and generated by a single element, and conversely every closed subring with these properties is of the form $L(\phi)$ for some $\phi$ in $C(K)$.

Thus saying that $\basic (K) \le \kappa$ is the same as saying that $C(K)$ is the sum of no more than $\kappa$ closed subrings containing the constants and generated by a single element. So the results above imply that the problem of deciding whether the Banach algebra $C(K)$ can be written as a sum of a certain size of `small' closed subrings is closely linked to $\cof ([w(K)]^{\aleph_0} , \subseteq)$ and PCF theory.

\section{Separable Metrizable Spaces}

The following simple lemma is used repeatedly and without further reference.
Let $\bfam$ be a basic family for a space $X$, and let $C$ be a $C$--embedded subspace (every continuous real valued function on $C$ can be extended over $X$). Then clearly $\bfam \restriction C = \{ \phi \restriction C : \phi \in \bfam\}$ is basic for $C$. Hence:
\begin{lemma}
Let $C$ be a $C$--embedded subspace of a space $X$ --- for example if $X$ is normal, and $C$ is closed --- then $\basic (X) \ge \basic (C)$.
\end{lemma} 

\begin{theorem}\label{sm}
Let $X$ be separable metrizable. Then either $\basic (X)$ is finite, which occurs if and only if $X$ is locally compact and finite dimensional, or $\basic (X) = \ctm$.
\end{theorem}

\begin{proof}
Let $X$ be separable metrizable. Four cases arise.

The first case is when $X$ is locally compact and finite dimensional. Then $\basic (X) \le 2 \dim (X)+1$, by the Main Theorem of \cite{FG1}.

In all remaining cases we show $\basic(X) \ge \ctm$, and so equals the continuum.

The second case is when $X$ is not locally compact. Then, as $X$ is first countable and normal, $X$ contains a closed copy of the metric fan, $F$ (defined below). So $\basic (X)\ge \basic (F) \ge \ctm$ by Proposition~\ref{fan} and Proposition~\ref{n_bdd}.

Case~3 is that $X$ is locally compact, infinite dimensional,  but contains no infinite dimensional compact subspaces. Then we can write $X$ as a union of open sets $(U_n)_n$ such that, for all $n$,  compact $\cl{U_n} \subset U_{n+1}$ and $\dim (U_n) < \dim (U_{n+1})$. Using the Countable Sum Theorem for dimension, we can extract compact subsets $C_n$ from the `gaps' $U_{n+1} \setminus \cl{U_n}$ such that $\dim C_n < \dim C_{n+1}$ for all $n$. Now we see that $C$, the disjoint union of the $C_n$'s is a closed subspace of $X$ satisfying the conditions of Proposition~\ref{incing}, so we indeed have, $\basic (X) \ge \basic (C) \ge \ctm$.

Finally, suppose $X$ is locally compact and contains an infinite dimensional compact subspace $K$. It suffices to show $\basic (K) \ge \ctm$, which is the content of Proposition~\ref{cpt_inf}.
\end{proof}

\paragraph{Independent Families}

In vector spaces one method of giving a lower bound for the size of a basis is to find large linearly independent sets. We apply the same approach to give lower bounds for $\basic(X)$. Note that if $V$ is a vector space, then $L \subseteq V$ is linearly independent if and only if its intersection with any subspace spanned by $n$ members of $V$ contains no more than $n$ elements. This leads us to the correct definition of `functional independence'.

Let $\mathcal{C}$ be a subset of $C(X)$. We say that $\mathcal{C}$ is (functionally) {\sl independent} if for all $n$, and any $\phi_1, \ldots ,\phi_n \in C(X)$ we have $|\mathcal{C} \cap \{ \sum_{i=1}^n g_i \circ \phi_i : g_1 , \ldots , g_n \in C(\mathbb{R})\}| \le n$. (We  omit the adjective `functionally' except when we need to differentiate from linear independence in the vector space sense.)

Further, we say $\mathcal{C}$ is {\sl weakly independent} if for all $n$, and any $\phi_1, \ldots ,\phi_n \in C(X)$ we have $|\mathcal{C} \cap \{ \sum_{i=1}^n g_i \circ \phi_i : g_1 , \ldots , g_n \in C(\mathbb{R})\}| < \mathfrak{c}$, and we say $\mathcal{C}$ is {\sl strongly independent} if for all $n$, and any $\phi \in C(X, \mathbb{R}^n)$ we have $|\mathcal{C} \cap \{ g \circ \phi : g \in C(\mathbb{R}^n)\}| \le n$.

Clearly `independent' implies `weakly independent'. Further, writing $\sum_{i=1}^n g_i \circ \phi_i$ as $g \circ \phi$ where $\phi (x_1, \ldots , x_n) = (\phi_1 (x_1) \ldots  , \phi_n (x_n))$ and $g(y_1, \ldots , y_n) = \sum_{i=1}^n g_i (y_i)$, we see that ``strongly independent' implies `independent'.

\begin{lemma}
If a space $X$ has a weakly independent family $\mathcal{C}$ of size $\ge \mathfrak{c}$, then $\basic (X) \ge \mathfrak{c}$.
\end{lemma}
\begin{proof}
Let $\bfam$ be a basic family for $X$. For each $f \in \mathcal{C}$, pick $\phi_1, \ldots , \phi_n$ from $\bfam$ so that $f = \sum_{i=1}^n g_i \circ \phi_i$. Then as $\mathcal{C}$ is weakly independent, the map taking $f$ in $\mathcal{C}$ to $\{ \phi_1 , \ldots , \phi_n\}$ in $\bigcup_{m \in \mathbb{N}} [\bfam]^m$ is $< \mathfrak{c}$--to--$1$. Since $|\mathcal{C}| \ge \mathfrak{c}$, it follows that $|\bfam| \ge \mathfrak{c}$ --- as required.
\end{proof}

To create large functionally independent families we will start from large linearly independent sets in the vector space $\mathbb{R}^n$ (with its usual inner product).

\begin{prop}\label{embs}
Fix a natural number $n$.
\begin{itemize}
\item[(a)] There is a Cantor set $C$ contained in the unit $(n-1)$--sphere of $\mathbb{R}^n$ such that for any distinct $x_1, \ldots, x_n$ in $C$, the $x_i$'s form a basis of $\mathbb{R}^n$.

\item[(b)] Let $J$ be a non--trivial closed bounded interval,  and $B$ a homeomorph of the $n$--cube, $J^n$.  There is a Cantor set $D$ contained in $C(B,J)$ such that  for any distinct $d_1, \ldots , d_n$ in $D$ the map $d=(d_1, \ldots , d_n) : B \to J^n$ is an embedding.
\end{itemize}
\end{prop}

\begin{proof}[Proof (of (a))]
Let $U=\{ (x_1, \ldots , x_n) \in (\mathbb{R}^n)^n : x_1 , \ldots , x_n$ are linearly independent$\}$. Then $U$ is open and dense in $(\mathbb{R}^n)^n$. One can further check that $U_K = \{ K \in \mathcal{K}(\mathbb{R}^n) : \text{ for all distinct } x_1, \ldots , x_n \in K \ (x_1, \ldots , x_n) \in U\}$ is comeagre in the space  $\mathcal{K}(\mathbb{R}^n)$ of compact subsets of $\mathbb{R}^n$ with the Hausdorff metric. (This is the Mycielski--Kuratowski technique, see 19.1 of \cite{Kechris}.) Since the set $P_K = \{ K \in \mathcal{K}({\mathbb{R}^n}) : K \text{ is perfect}\}$ is also comeagre in the Polish space  $\mathcal{K}(\mathbb{R}^n)$, and perfect compact metric spaces contain Cantor sets, we can indeed pick a  Cantor set $C \subseteq \mathbb{R}^n$ such that for any distinct $x_1, \ldots, x_n$ in $C$, the $x_i$'s are linearly independent, and hence form a basis. Mapping each $x$ in $C$ to $x/ \|x\|$ we see we can assume $C$ is contained in the unit $(n-1)$--sphere.
\end{proof}

\begin{proof}[Proof (of (b))] First note that if (b) holds for one choice of $J$ and $B$, then it holds for all. We will  use the interval $J=[-1,+1]$, and the closed $n$--ball, $B^{(n)}$. Also note that we work in the inner product space $\mathbb{R}^n$.

Fix a Cantor set $C$ in the unit sphere of $\mathbb{R}^n$ as in part~(a). Let $\hat{C} = \{ \hat{c} : c \in C\}$ where $\hat{c}$ is the linear functional on $\mathbb{R}^n$ dual to $c$, namely $\hat{c} (x) = \langle c , x\rangle$. Then, by duality, $\hat{C}$ is a Cantor set in $\mathbb{R}^* \subseteq C(\mathbb{R}^n, \mathbb{R})$, and any $n$--many distinct elements of $\hat{C}$ are linearly independent.

Let $D= \{ \hat{c} \restriction    B^{(n)} : c \in C\}$. Then $D$ is a family of continuous functions mapping $B^{(n)}$ to $[-1,+1]$, with the required properties.
\end{proof}

\paragraph{Compact Case, Fixed $n$}

\begin{prop}\label{cpt_dim>n}
Fix $K$  a compact  space of dimension  $>n \ge 2$. 

Then there is a Cantor set $C \subseteq C(K,I)$ such that for all $\phi \in C(K,I^n)$ we have $|C \cap \{ g \circ \phi : g \in C(I^n,I)\}| \le n$.
\end{prop}

\begin{proof} Recall (see \cite{Alex}, for example) that a normal space, $X$, has dimension $\le n$ if and only if every continuous map from a closed subspace into the $n$--sphere (which is homeomorphic to the boundary of the $(n+1)$--cube) has a continuous extension over $X$.
Hence, as $\dim K >n$,  there is a map $p : K \to I^{n+1}$ and closed subspace $A$, such that $p \restriction   A : A \to \partial I^{n+1}$ can not be continuously extended (over $K$ into $\partial I^{n+1}$). We may suppose that $A = p^{-1} \partial I^{n+1}$.

By Proposition~\ref{embs}~(b) there is a Cantor set $D$ contained in $C(I^{n+1},I)$ such that for any distinct $d_1, \ldots , d_{n+1} \in D$ the map $d=(d_1, \ldots , d_{n+1}):I^{n+1} \to I^{n+1}$ is an embedding. 
For distinct $d_1, \ldots , d_{n+1} \in D$, and embedding $d=(d_1, \ldots , d_{n+1})$ define $f_d = d \circ p$. Note that $f_d \ne f_{d'}$ if $d \ne d'$. Let $C= \{ f_d : d \in D\}$. This is a Cantor set in $C(K,I^{n+1})$.

Suppose, for a contradiction, for some $\phi \in C(K,I^n)$, there were $(n+1)$ distinct elements $f_1, \ldots , f_{n+1}$ in $C \cap \{ g \circ \phi : g \in C(I^n,I)\}$. So,  for  $i=1, \ldots , n+1$, we have  $f_i = d_i \circ p$ for some (distinct) $d_i \in D$, and $f_i = g_i \circ \phi$  for some $g_i \in C(I^n,I)$.

Let $d=(d_1, \ldots , d_{n+1})$, and $g=(g_1, \ldots , g_{n+1})$. So $p \circ d = g \circ \phi$. Since $d$ is an embedding, we have $p = h \circ \phi$ where  $h= (d^{-1} \circ g)$ is in $C(I^n,I^{n+1})$.

Let $A' = h^{-1} \partial I^{n+1}$. Note that $\phi^{-1} A' = p^{-1} \partial I^{n+1} =A$, so $\phi$ maps $A$ inside $A'$. Since $K' = \phi (K)$ is contained in $I^n$ it has dimension $\le n$. Hence the map $h \restriction    A' : A' \to \partial I^{n+1}$ has a continuous extension $h' : K' \to \partial I^{n+1}$.

But now $p \restriction   A : A \to \partial I^{n+1}$ has a continuous extension over $K$ into $\partial I^{n+1}$ --- namely $h' \circ \phi$ --- contradiction!
\end{proof}

\paragraph{Locally Compact, All Compact Subspaces Small}

\begin{prop}\label{incing}
Let $(C_n)_n$ be a sequence of compact spaces  such that each $C_n$ has finite dimension $>n$. Let $X = \bigoplus_n C_n$, and $\gamma X$ be a  compactification of $X$.

Then there is a Cantor set $C$ contained in $C(\gamma X,I) \subseteq C(X)$ such that
$C$ is strongly independent for $C(X)$ (and hence for $C(\gamma X)$).

Hence $\basic (X) \ge \ctm$ and $\basic (\gamma X) \ge \ctm$
\end{prop}

\begin{proof}
For each $n \ge 2$, fix the Cantor set, $E_n$, guaranteed by Proposition~\ref{cpt_dim>n} for the $>n$ dimensional space $C_n$, and fix a homeomorphism $h_n$ from the standard Cantor set $\mathbf{C}$ to $E_n$. Let $C= \{ f_c : c \in \mathbf{C}\}$ where $f_c$ is constantly equal to zero on $C_1$ and on the remainder $\gamma X \setminus X$, and equals $h_n(c)/n$ on $C_n$. Note that each $f_c$ is continuous, and so $C$ is a Cantor set in $C(\gamma X,I)$.

Take any $n \ge 2$ and $\phi \in C(X, \mathbb{R}^n)$. Considering the restrictions of $\phi$ and elements of $C$ to $C_n$, it is immediate from the properties of $E_n$, that $|C \cap \{ g \circ \phi : g \in C(\mathbb{R}^n) \}| \le n$. Thus $C$ is strongly independent.
\end{proof}

\paragraph{Compact, Infinite Dimensional}

\begin{prop}\label{cpt_inf} Let $K$ be compact and infinite dimensional. Then there is a Cantor set $C$ contained in $C(K,I)$ which is strongly independent.

Hence, $\basic (K) \ge \ctm$.
\end{prop}

\begin{proof}
We show an appropriate, strongly independent, Cantor set $C$ exists. 
Dowker has shown \cite{Dowker} that if $X$ is a normal space and $M$ is a closed subspace with $\dim \le n$ then $\dim X \le n$ if and only if $\dim F \le n$ for all closed subsets of $X$ disjoint from $M$. In particular: \ ($\ast$) \ if $M$ contains a single point, $x$, then $\dim X >n$ if and only if $\dim F >n$ for some closed subset $F$ of $X \setminus \{x\}$.
For each point $x$ in $K$ pick a closed neighborhood of minimal dimension, $B_x$. By compactness, for some $x$, $B_x$ is infinite dimensional, and so all  neighborhoods of $x$ are infinite dimensional. Let $K_1=K$. Apply ($\ast$) to get a compact subset $C_1$ of $K_1$ not containing $x$ with $\dim C_1 >1$. Pick a closed neighborhood $K_2$ of $x$ disjoint from $C_1$. Apply ($\ast$) to get a compact subset $C_2$ of $K_2$ not containing $x$ with $\dim C_2 > \max (2,\dim C_1)$. Inductively, we get a 
pairwise disjoint collection, $\{ C_n : n \in \mathbb{N}\}$, of compact subsets of $K$ which are either (i) of strictly increasing (finite) dimensions, or (ii) all infinite dimensional. Let $K'$ be the closed subspace $\overline{\bigoplus_n C_n}$.

In the first case we apply Proposition~\ref{incing} to $K'$  to get a strongly independent Cantor set in $C(K')$ -- and hence in $C(K)$ -- as required.

In the second case, by Proposition~\ref{cpt_dim>n}, for each $n$ there is a Cantor set $E_n \subseteq C(C_n,I)$ such that for all $\phi \in C(C_n,I^n)$ we have $|E_n \cap \{g \circ \phi : g \in C(I^n,I)\}|$ finite. Fix homeomorphisms $h_n$ between the standard Cantor set $\mathbf{C}$ and $E_n$.

Define, for $c \in \mathbf{C}$, a map $f_c : K' \to I$ by: $f_c$ is identically zero on $K' \setminus \bigoplus_n C_n$ and $f_c(x') = (1/n) h_n(c) (x')$ if $x' \in C_n$. Then the $f_c$'s are continuous, can be continuously extended over $K$, and so form a Cantor set $C$ in $C(K,I)$. Further, if $\phi \in C(K,I)$ and $f_1, \ldots, f_{n+1} \in C$, then the $f_i$'s are not all in $\{ g \circ \phi : g \in C(I^n,I)\}$, because $f_1 \restriction E_n , \ldots , f_{n+1} \restriction    E_n$ are not all in $\{ g \circ (\phi \restriction   E_n) : g \in C(I^n,I)\}$, by choice of $E_n$. 

Thus the Cantor set $C$ is strongly independent as required.
\end{proof}

\paragraph{The Non Locally Compact Case} 
Let $F$ be the \emph{metric fan} where $F=
(\mathbb{N} \times \mathbb{N}) \cup \{\ast\}$, points in $\mathbb{N}
\times \mathbb{N}$ are isolated and basic neighborhoods of $\ast$ are
$B(\ast, n) = ([n,\infty) \times \mathbb{N}) \cup \{\ast\}$. Then a
  separable metric space is not locally compact if and only if it
  contains a closed copy of the metric fan. Thus if $\basic (F) =
  \mathbf{c}$ then $\basic(X)=\mathbf{c}$ for every separable metric
  space $X$ which is not locally compact.

We first reduce the calculation of $\basic (F)$ to that of
$\basic (\mathbb{N}, [-1,+1])$. Here we say that a family $\hat{\bfam}
\subseteq C(\mathbb{N},[-1.+1])$ is `basic for $\mathbb{N}$ into $[-1,1]$' if 
$\forall \hat{f} \in
C(\mathbb{N},[-1,+1])$  there are $\hat{ \phi}_1, \ldots , \hat{ \phi}_n \in
\hat{\bfam}$, and  $\hat{g}_1, \ldots , \hat{g}_n \in C(\mathbb{R})$ such that $\hat{f} = \sum_{i=1}^n \hat{g}_i
\circ \hat{ \phi}_i$, and define $
 \basic (\mathbb{N},[-1,+1]) = \min \{ |\hat{\bfam}| : \hat{\bfam}$ is basic for $\mathbb{N}$ into $[-1,1]\}$.

\begin{prop}\label{fan}
$\basic (F) \ge \basic (\mathbb{N},[-1,+1])$.
\end{prop}
\begin{proof}
Let $\bfam$ be basic for $F$. We will show that there is a
$\hat{\bfam}$ with $|\hat{\bfam}|=|\bfam|$ such that $\hat{\bfam}$ is
basic for $\mathbb{N}$ into $[-1,+1]$. 

For each $ \phi \in \bfam$ and $n$ such that $ \phi$ maps $\{n\}\times
\mathbb{N}$ into $[-1,+1]$, define $\widehat{ \phi_n}$ in
$C(\mathbb{N},[-1,+1])$ by $\widehat{ \phi_n} (m) =  \phi(n,m)$. Let
$\widehat{\bfam_n} = \{ \widehat{ \phi_n} :  \phi \in \bfam\}$ and
$\widehat{\bfam}= \bigcup_n \widehat{\bfam_n}$. Note that
$|\widehat{\bfam}|=|\bfam|$. 

Take any $\hat{f} \in C(\mathbb{N},[-1,+1])$. Define $f : F \to
[-1,+1]$ by $f(\ast)=0$ and $f(n,m)= \hat{f}(m)/n$. Note $f$ is
continuous. So there are $ \phi_1 , \ldots ,  \phi_n$ in $\bfam$ and 
$g_1, \ldots , g_n$ in $C(\mathbb{R})$ such that $f = \sum_i g_i \circ  \phi_i$.

By continuity of $ \phi_1, \ldots ,  \phi_n$ at $\ast$ there is an $N$
such that each $ \phi_i$  maps $\{N\} \times \mathbb{N}$
into a closed bounded interval, say $I_i$. Fix  homeomorphisms $h_i$ of
$\mathbb{R}$ with itself carrying $I_i$ to $[-1,+1]$. Now we see that,
replacing $g_i$ with $g_i \circ h_i^{-1}$ and $ \phi_i$ with $h_i \circ
 \phi_i$, we can assume that the $ \phi_i$ all map into $[-1,+1]$.

Thus $\widehat{ \phi_1} = \widehat{ ( \phi_{1})_N},
  \ldots , \widehat{ \phi_n} = \widehat{( \phi_n)_N}$ are in
  $\widehat{\bfam}_N \subseteq \widehat{\bfam}$. Further, as
  $\hat{f}(m)/N = f(N,m) = \sum_{i=1}^n g_i( \phi_i (N,m)) = \sum_i
  g_i(\widehat{ \phi_i}(m))$, we have that $\hat{f} = \sum_{i=1}^n
  \widehat{g_i} \circ \widehat{ \phi_i}$ where $\widehat{g_i}=N. g_i$
  --- as required. 
\end{proof}

\begin{prop}\label{n_bdd}
 There is a Cantor set $C$ contained in $C(\mathbb{N},[-1,+1])$ such that $|C \cap \{\sum_{i=1}^n g_i \circ \phi_i : g_1, \ldots , g_n \in C(\R)\}| \le {\aleph_0}$ for all  $\phi_1, \ldots , \phi_n$ from  $C(\mathbb{N},[-1,+1])$.

Thus $C$ is `weakly independent' in the sense appropriate for $C(\mathbb{N},[-1,+1])$, and so $\basic (\mathbb{N},[-1,+1]) = \mathbf{c}$.
\end{prop}
\begin{proof}
 Define $C= \{ f \in C(\mathbb{N},[-1,+1]) : 
f(\mathbb{N}) = \{-1,+1\}\}$.  Then $C$ is a Cantor set, and we
will prove that, for each $n$, and finite $\bfam' \subseteq C(\mathbb{N}, [-1,+1])$ we have  $|C \cap L(\bfam')| = {\aleph_0}$. 

Fix $n \ge 1$. Fix $\phi \in C(\mathbb{N},[-1,+1]^n)$. As in the argument that `strongly independent' implies `independent' to prove the claim it suffices to show that there are only
countably many $f \in C$ representable as $g \circ \phi$ for some $g \in C([-1,+1]^n,[-1,+1])$. 

Let $K=\cl{\phi(\mathbb{N})}$ --- a compact subset of
$[-1,+1]^n$. A composition $g \circ \phi : \mathbb{N} \to [-1,+1]$ is determined by the values
of $g$ on $\phi (\mathbb{N})$, and so definitely determined by its values on $K$.

If $g \circ \phi$ is in $C$, then, by continuity, $g
\restriction K$
maps $K$ onto $\{-1,+1\}$. Thus $K$ is partitioned into two non--empty
clopen pieces, one of which is mapped by $g$
to $-1$, and the other to $+1$. But a compact metric space only has
countably many clopen subsets. So there are only a countable number of
possibilities for $g$ on $K$, and only
countably many $f \in C$ representable as $g \circ \phi$ --- as claimed.
\end{proof}

\begin{cor} Let $X$ be finite dimensional, locally compact, \textbf{not compact}, separable metrizable. Then:

(1) \ there is a basic family $\bfam \subseteq C(X)$ such that $\bfam$ is finite, but

(2) \ there is no basic$^*$ family $\bfam^*$ consisting of \textbf{bounded} functions such that $|\bfam^*| < \mathbf{c}$.
\end{cor}

\begin{proof} The first claim is just the Main Theorem of \cite{FG1}. For the second part, first note that since $\N$ can be embedded as a closed subspace of
$X$, it is sufficient to show that (2) holds for $\N$.
Suppose, for contradiction, there exists a basic family $\bfam^*$
for $\N$ consisting of bounded function whose cardinality is $<\mathbf{c}$.

Write $\bfam^*=\bigcup_{n\in\N}\bfam_n$ where
$\bfam_n=\{\phi:-n\leq\phi(a)\leq n,\text{ for each }n\in\N\}$. Then
$C^*(\N)=\bigcup_{n\in\N}L(\bfam_n)$. Let $\mathcal{F}= \{ f \in
C(\mathbb{N},[-1,+1]) : f(\mathbb{N}) = \{-1,+1\}\}$ as in the proof
of Proposition~\ref{n_bdd}. There exists an $m_0$ such that
$|\mathcal{F}\cap L(\bfam_{m_0})|=\mathbf{c}$.  But the argument in the proof
of Proposition~\ref{n_bdd} shows $L(\bfam_{m_0})\leq |\bfam^*|<\mathbf{c}$
which is the desired contradiction.
\end{proof}

\section{Compact Spaces}

\begin{prop}\label{cpt_fd}
Suppose $K$ is compact and finite dimensional. Then $\basic (K) \le \cof
([w(K)]^{\aleph_0} , \subseteq)$.
\end{prop}

\begin{proof}
Let $K$ be compact of dimension $n$. Then there is a directed set
$(\Lambda, \le)$ where $|\Lambda| = w(K)$, compact metric $K_\lambda$
with $\dim K_\lambda \le n$, and for all $\lambda \ge \mu$ a continuous
map $f_{\lambda, \mu}$ such that $K = \plim \{K_\lambda : \lambda \in
\Lambda\} = \{ \langle
x_\lambda \rangle \in \prod_{\lambda} K_\lambda: \lambda \ge \mu
\implies f_{\lambda , \mu} (x_\lambda) =x_\mu\}$.

Let $\mathcal{C}$ be cofinal in $([w(K)]^{\aleph_0} , \subseteq)$. We may
suppose that each $C$ in $\mathcal{C}$ is directed.
For each $C \in \mathcal{C}$, $K_C = \plim \{K_\lambda : \lambda \in
C\}$ is compact, metric of dimension $\le n$. So $K_C$ has a basic
family $\bfam_C'$ of size $2n+1$. Define $p_C = \pi_C \restriction    \plim
\{K_\lambda : \lambda \in \Lambda\}$. 
Define $\bfam_C = \{ \phi' \circ p_C
:  \phi' \in \bfam_C'\}$. and $\bfam = \bigcup_{C \in \mathcal{C}}
\bfam_C$.
Then $|\bfam| = |\mathcal{C}|$. We show that $\bfam$ is basic -- as
required. 

To this end, take any $f \in C(K)$.
Extend $f: \plim \{x_\lambda : \lambda \in \Lambda\} \to \R$ to
continuous $\hat{f} : \prod_{\lambda \in \Lambda} K_\lambda \to
\R$. Then there is a countable $\Lambda_0 \subseteq \Lambda$ and
continuous $g_0 : \prod_{\lambda \in \Lambda_0} K_\lambda \to \R$ such
that $\hat{f} = g_0 \circ \pi_{\Lambda_0}$. 
Pick $C \in \mathcal{C}$
such that $C \supseteq \Lambda_0$. Note that as $C$ is directed,
$\{\langle x_\lambda \rangle_{\lambda \in C} : \lambda \ge \mu
\implies f_{\lambda , \mu} (x_\lambda)=x_\mu\} = \plim \{K_\lambda :
\lambda \in C\}$, and $\pi_C$ maps $\plim \{K_\lambda : \lambda \in
\Lambda\}$ to $\plim \{K_\lambda : \lambda \in C\}$. 

We can write $\hat{f} = \hat{g} \circ \pi_C$ where $\hat{g} = g_0
\circ \pi^C_{\Lambda_0}$ is a continuous map $\prod_{\lambda \in C}
K_\lambda$ into $\R$. Thus $f = \hat{f} \restriction    \plim \{K_\lambda : \lambda \in \Lambda\} = g
\circ p_C$ where $p_C = \pi_C 
\restriction    \plim \{ K_\lambda : \lambda \in
\Lambda\}$ and $g = \hat{g} 
\restriction    \plim \{ K_\lambda : \lambda \in C\}$.

Now we see that $g = \sum_{i=1}^{2n+1} g_i \circ  \phi_I'$ where
$ \phi_C' \in \bfam_C'$ and $g_i \in C(\R)$. Thus 
\[ f = g \circ p_C = \sum_{i=1}^{2n+1} g_i \circ \left(  \phi_i' \circ
\pi_C\right) = \sum_{i=1}^{2n+1} g_i \circ  \phi_i, \]
where $ \phi_1, \ldots ,  \phi_{2n+1}$ are in $\bfam_C \subseteq \bfam$
and $g_1, \ldots , g_{2n+1}$ are in $C(\R)$.
\end{proof}

Suppose $K$ is compact and $w(K)$ has uncountable cofinality. Then recalling that $\cof ([w(K)]^{\aleph_0}, \subseteq) =\kappa$ in this case, from Propositions~\ref{cpt_inf} and~\ref{cpt_fd} we deduce:

\begin{theorem}\label{cpt_dichot1}
If $K$ is compact and its weight has uncountable cofinality, then 
\begin{itemize}
\item[either] $K$ is finite dimensional, $\basic (K) = \cof
([w(K)]^{{\aleph_0}}) = w(K)$, and $\basic (K) < w(K)^{{\aleph_0}}$ if and
only if $w(K) < \mathfrak{c}$,
\item[or] $K$ is infinite dimensional, $\basic (K) = |C(K)| = w(K)^{\aleph_0}$, and $w(K) < \basic (K)$ if and
only if $w(K) < \ctm$,
\end{itemize}
\end{theorem}

Thus, considering only compact spaces $K$ whose weight has uncountable cofinality, the statements: `there is a space with $\basic (K) < w(X)^{\aleph_0}$', `there is a space with $w(K) < \basic (K)$', and `the continuum hypothesis fails', are all equivalent. Further, `there is a space with $w(K) < \basic (K) < w(K)^{\aleph_0}$' is false.

Call a space $X$ `nice' if it contains a discrete subset $D$ with $|D|=w(X)$.
Note that there are many examples of compact `nice' spaces, for example: $2^\kappa$, $I^n \times 2^\kappa$ and $I^\kappa$ are compact, `nice' and span the dimensions.

\begin{prop}\label{cpt_nice} If $K$ is compact and `nice', then $\basic (K)\geq
\cof([w(K)]^{{\aleph_0}},\subseteq)$.
\end{prop}

\begin{proof} Let $D$ be discrete in $K$ with $w(K)=|D|$. Let $K'=\overline{D}$, and $K'_c = K'\setminus D$. Since $w(K')=w(K)$ and $\basic (K) \ge \basic (K')$ it suffices to show $\basic(K') \ge \cof([w(K')]^{{\aleph_0}},\subseteq)$.

Note that $D$ is open in $K'$, so $K'_c$ is compact. Take any function $f\in C(K', \R^n)$. Then  $f(K'_c)$ is a
compact subset of $\R^n$, so it is a $G_{\delta}$ subset, and we can write $f(K'_c)$ as $\bigcap_{n\in\N}U_n$,
where $U_n$ is open set in $\R^n$ for each $n$. As $K'$ is compact, each $K'\setminus f^{-1}(U_n)$ is closed and discrete, and hence
finite. So we can define a countable subset of $D$ for
each $f$ by $C_f=\bigcup_{n\in\N} K'\setminus f^{-1}(U_n)$.

Now suppose $\bfam\subseteq C(K')$ with 
$|\bfam|<\cof(|w(K')|^{\aleph_0},\subseteq)$. We show $\bfam$ is not
a basic family.

Given $\phi_1,\phi_2,\cdots, \phi_n$ from $\bfam$, let $\hat{\phi}=(\phi_1, \ldots , \phi_n) : K' \to \R^n$, and $C(\phi_1, \ldots , \phi_n) = C_{\hat{\phi}}$.
 Let $\mathcal{C}=\{C(\phi_1, \ldots , \phi_n): \phi_1, \ldots , \phi_n \in \bfam\}$.
Since $|\bfam|<\cof([w(K')]^{\aleph_0},\subseteq)$, the
collection $\mathcal{C}$ is not cofinal in $[D]^{\aleph_0}$. Therefore
there exists a countably infinite subset $C$ of $D$ such that for any
$\phi_1,\cdots, \phi_n$, $C$ is not a subset of $C(\phi_!, \ldots , \phi_n)$. 

 Take any $\phi_1,\ldots, \phi_n$ in $\bfam$. Pick $x$ in $C$ but not $C(\phi_1, \ldots , \phi_n)$. By definition of $C(\phi_1, \ldots , \phi_n)$ there exists $x'\in K'_c$ such that
$\hat{\phi}(x)=\hat{\phi}(x')$. Then for any $g_1, \ldots , g_n$ from $C(\R)$,
$\sum_{i=1}^n g_i \circ \phi_i$ takes the same value at a point in $C$ and at a point in $K'_c$. 

But now we see that if we enumerate $C= \{x_1, x_2, \ldots \}$, and define $h$ by $h(x_n)=1/n$ and $h$ is identically zero outside $C$, then $h$ is continuous and $h(C)$ is disjoint from $h(K'_c)$. Thus $h$ can not be represented by any finite
collection of $\bfam$, and so $\bfam$ is not  basic.
\end{proof}

From the identity $w(K)^{\aleph_0} = \cof ([w(K)]^{\aleph_0}, \subseteq) \times \mathfrak{c}$ and Propositions~\ref{cpt_inf}, \ref{cpt_fd} and~\ref{cpt_nice} we conclude:

\begin{theorem}\label{cpt_dichot2}
If $K$ is compact and `nice' then:
\begin{itemize}
\item[either] $K$ is finite dimensional and $\basic (K) = \cof ([w(K)]^{\aleph_0} , \subseteq)$,
\item[or] $K$ is infinite dimensional and $\basic (K) = |C(K)|=w(K)^{\aleph_0}$.
\end{itemize}
\end{theorem}

Thus considering only `nice', finite dimensional, compact spaces $K$ whose weight has countable cofinality (for example, $K=2^{\aleph_\omega}$), it is always true that $w(K) < \basic (K)$, and it is consistent and independent (depending on the value of the continuum, $\mathfrak{c}$) whether $\basic (K) < w(K)^{\aleph_0} = |C(K)|$.

\section{Open Questions}

The most immediate question is whether the restriction to `nice' compacta in Proposition~\ref{cpt_nice} is necessary.
\begin{ques}
Is it true that  $\basic (K)\geq
\cof([w(K)]^{{\aleph_0}},\subseteq)$ for all compact spaces $K$?
\end{ques}

The proofs of the results for compact spaces clearly rely on facts and techniques that only apply to compact spaces. But it seems possible that the results could be extended to larger classes of spaces.
\begin{ques}
Do the results for $\basic (K)$ for compact $K$ hold for (1) locally compact, Lindelof spaces or even (2) all Lindelof spaces? 
\end{ques} 

In a different direction, what about discrete spaces?
\begin{ques}
Is $\basic (D({\aleph_0}_1)) = \aleph_1$? $=2^{{\aleph_0}}$?
\end{ques}

\end{document}